\documentclass[12pt]{article}
\begin{document}
\begin{center}
{\large {\bf FINITE -- DIMENSIONAL
REPRESENTATIONS\\
OF  $U_q[gl(2/1)]$ 
IN A BASIS OF $U_q[gl(2)\oplus gl(1)]$}}\\[1cm] 
{\bf Nguyen Anh Ky}~$^{1,2}$ ~ and ~ {\bf Nguyen Thi Hong Van}~$^1$\\[5mm] 
$^1$ Institute of Physics\\ P.O. Box 429, Bo Ho, Hanoi 10 000, Vietnam\\[4mm]
$^2$ The Abdus Salam International Centre for 
Theoretical Physics (ICTP),  I-34100 Trieste, Italy

\end{center}
 \vspace*{5mm}	
\begin{abstract}

    The quantum superalgebra $U_q[gl(2/1)]$ is given as both  a Drinfel'd--Jimbo 
deformation of $U[gl(2/1)]$  and a Hopf superalgebra. Finite--dimensional 
representations of this quantum superalgebra are constructed and investigated in a 
basis of its even subalgebra $U_q[gl(2)\oplus gl(1)]$.  The present method for 
constructing representations of a quantum superalgebra combines previously suggested 
ones for the cases of superalgebras and quantum superalgebras, and, therefore, has 
an advantage in comparison with the latter.  
\end{abstract} 
\vspace*{4mm}

\underline{PACS numbers} : ~ 02.20Tw, 11.30Pb.\\
 
\underline{MSC--class.} : 81R50; 17A70.\\[1cm]
{\bf I. INTRODUCTION}

   Emerged about twenty years ago \cite{frt} -- \cite{woro} from the study on the 
quantum inverse scattering method  and Yang--Baxter equations \cite{collect}, 
quantum  groups  (QG's) readily became one of the most interesting concepts in 
physics and mathematics in the last two decades. 
For a short time QG's and their representations have been investigated in 
details in both the physical and the mathematical aspects and have  found 
various applications in physics \cite{2dim} -- \cite{majid}.

  One of the approaches to QG's is the Drinfel'd--Jimbo (DJ), or "quantum", 
deformation of universal enveloping algebras \cite{drin,jim}. This kind of 
deformation depends on one or more parameters which could be generic complex 
numbers or roots of unity.  The defined in this way QG's appear to be Hopf algebras 
which are typically noncommutative and noncocommutative \cite{drin}. The latter, 
in turn,  can be used for introducing and studying QG's.  Hopf algebra structures of 
QG's are shown to be an efficient tool for investigating QG's as the whole and their 
representations in particular.  Moreover, these investigations can be extended to 
quantum supergroups (QSG's), a notion combining QG's with supersymmetry 
\cite{kure} -- \cite{khoro}. For their generality and importance QSG's, e.g.,  
$U_q[gl(m/n)]$,  which are deformations of  universal enveloping algebras 
$U[gl(m/n)]$ of  superalgebras $gl(m/n)$ are a subject of research interest 
\cite{manin3} -- \cite{pq01}. Representations of these QSG's called also 
quantum superalgebras (QSA's)  are explicitely known in a number of cases but 
their conctructions are sometimes complicated with heavy calculations, especially 
for higher rank cases. We suggested in  \cite{q94} a method (procedure) for 
constructing and investigating finite--dimensional representations of  a QSA. 
This method is very efficient for one-parametric deformations $U_q[gl(m/n)]$, at 
least with $m$ and $ n$ not very high \cite{q94} -- \cite{ic96},  and it can be also  
applicable to the two-parametric case \cite{pq96,pq00,pq01}. In general, the method 
proposed is good, however,  as in the classical, i.e., non-deformed, case 
\cite{class89} -- \cite{osp32}, 
its practical realisation is not always convenient because the calculation process 
based on using (deformed) commutation relations between generators (without 
using their Hopf algebra structures) is cumbersome in some stages.  Besides that, 
the latter method does not give us an easy way to get an explicit description of 
representations of a QSA in a basis of its even subalgebra (for example, it does 
not express the so-called induced basis of a QSA in terms of a basis of the even 
subalgebra and, therefore, we do not have matrix elements in the induced basis 
as we could do in the classical case \cite{class89,sl13}). Such a description may 
be physically necessary as  in it both the origin and the structure of multiplets can  
be seen explicitly. 

   Exploiting the Hopf algebra structure of quantum superalgebras $U_q[gl(m/n)]$ 
we can investigate in a transparent and consistent way their module structure and 
representations. Taking a demonstration on $U_q[gl(2/1)]$ (which can be applied 
to physical problems such as those of strongly correlated electron systems 
\cite{correl1} -- \cite{correl3} ) we construct its induced module and find all its 
finite--dimensional irreducible representations in a basis of the even subalgebra 
$U_q[gl(2/1)_0]\equiv U_q[gl(2)\oplus gl(1)]$. The results obtained are hopefully 
useful for the above--mentioned  applications.  The present method combines the 
advantage of the previously suggested methods for the classical case 
\cite{class89,sl13} and the quantum deformation case \cite{q94,ic96} and may 
be  more convenient in some practical (calculation and application) aspects. This 
paper is organized as follows.  

   The quantum superalgebra $U_q[gl(2/1)]$ as a DJ deformation of $U[gl(2/1)]$ and 
as a Hopf superalgebra is given in Section II where its induced representations are also 
considered. Finite--dimensional representations of this quantum  superalgebra in a basis 
of its even subalgebra $U_q[gl(2)\oplus  gl(1)]$ (or  simply, a $U_q[gl(2)\oplus  
gl(1)]$ -- basis) are constructed in Section III and classified in Section IV. Finally, 
section V is devoted to some discussions and the conclusion. \\[1cm]  
{\bf II.  $U_q[gl(2/1)]$  AND  ITS  INDUCED  REPRESENTATIONS} 

   The quantum superalgebra $U_q[gl(2/1)]$ can be completely defined 
through the Weyl--Chevalley generators $E_{12}$,  $E_{21}$, $E_{23}$, 
$E_{32}$ and $E_{ii}$, $i=1,2,3$, which satisfy the following defining relations 
\cite{q94,ic96} :
\begin{tabbing} 
\=1234567891234567891\=$[E_{ii},E_{jj}]$1234\= =12\= 01234 
\=$[E_{ii},E_{j,j+1}]$\= =  
\=$(\delta_{ij}-\delta_{i,j+1})E_{j,j+1}$123\=\kill 

~~~~a)  {\it the super-commutation relations} ($1\leq i,i+1,j,j+1\leq 3$):\\[2mm] 
\>\>$[E_{ii},E_{jj}]$\> = \>0,\>\>\>\>(2.1a)\\[1mm] 
\>\>$[E_{ii},E_{j,j+1}]$\>=\>$(\delta_{ij}-\delta_{i,j+1})E_{j,j+1}$,\>\>\>\>(2.1b)\\[1mm] 
\>\>$[E_{ii},E_{j+1,j}]$\>=\>$(\delta_{i,j+1}-\delta_{ij})E_{j+1,j}$,\>\>\>\>(2.1c)\\[1mm] 
\>\>$[E_{12},E_{21}]$\>=\>$[H_1]$,\>\>\>\>(2.1d)\\[1mm] 
\>\>$\{E_{23},E_{32}\}$\>=\>$[H_2]$,\>\>\>\>(2.1e)\\[1mm] 
\>\>$H_i$\>=\>$(E_{ii}-{d_{i+1}\over 
d_{i}}E_{i+1,i+1}),$\>\>\>\>(2.1g)\\[4mm] 
where $d_{1}=d_{2}=-d_{3}=1$, and\\[4mm] 
~~~~b) {\it the Serre relations}:\\[2mm] 
\>\>~~~~~~$E_{23}^{2}$\>=\>~~~~$E_{32}^{2}$\>~~~~~~=~~~~~~0,\>\>\>(2.2a)\\[2mm] 
\>\>~~$[E_{12},E_{13}]_q$\>=\>$[E_{21},E_{31}]_q$\>~~~~~~~~=~~~~0, 
                                        \>\>\>(2.2b)
\end{tabbing}
with $E_{13}$ and $E_{31}$,  
$$E_{13}:=[E_{12},E_{23}]_{q^{-1}},\quad E_{31}:=-[E_{21},E_{32}]_{q^{-1}},
\eqno(2.3)$$
defined as new generators, where the notation
$$[A,B]_r:=AB-rBA$$
is used. The newly defined generators  are odd and have vanishing squares. The 
generators $E_{ij}$, $i,j=1,2,3$, including the new ones, are quantum deformation 
analogues ($q$--analogues) of the Weyl generators $e_{ij}$ of the classical 
superalgebra $gl(2/1)$ whose universal enveloping algebra $U[gl(2/1)]$ is a 
classical limit of $U_q[gl(2/1)]$  at $q\rightarrow 1$. The defined in this way 
quantum superalgebra $U_q[gl(2/1)]\equiv U_q$ is a Hopf superalgebra 
endowed with the following additional maps:
 
1) {\it coproduct}  
$\Delta$: ~~ $U_q $ $\rightarrow$ $U_q\otimes U_q$, 
$$\Delta (1)=1\otimes 1,$$ 
$$\Delta (E_{ii})=E_{ii}\otimes 1 + 1\otimes E_{ii},$$
$$\Delta (E_{12})=E_{12}\otimes q^{H_1}+ 1\otimes E_{12},$$
$$\Delta (E_{21})= E_{21}\otimes 1+ q^{-H_1}\otimes E_{21},$$
$$\Delta (E_{23})=E_{23}\otimes q^{H_2}+1\otimes E_{23},$$
$$\Delta (E_{32})= E_{32}\otimes 1+q^{-H_2}\otimes E_{32}.\eqno (2.4)$$

2) {\it antipode}  
$S$: ~~ $U_q$ $\rightarrow$ $U_q$,
$$S(1)=1,$$
$$S(E_{ii})=-E_{ii},$$
$$S(E_{12})=-E_{12}q^{-H_1},$$
$$S(E_{21})=- q^{H_1}E_{21},$$
$$S(E_{23})=-E_{23}q^{-H_2},$$
$$S(E_{32})=-q^{H_2}E_{32},\eqno (2.5)$$

3) {\it counit}  
$\varepsilon$: ~~ $U_q$ $\rightarrow$ ${\bf C}$,
$$\varepsilon (1)=1,$$
$$\varepsilon (E_{ii})=\varepsilon (E_{12})=\varepsilon (E_{21})=\varepsilon(E_{23}) 
=\varepsilon(E_{32})=0.\eqno (2.6)$$  
These maps are either homomorphisms ({\it coproduct} and {\it counit}) or 
an anti-homomorphism ({\it antipode}) and are consistent with the defining 
relations (2.1) -- (2.3). These relations are quantum deformations, or  
$q$-deformations, of the ordinary (super-) commutation relations and they 
can be obtained from the latter by replacing the classical adjoint operation 
$ad$ with the quantum deformation one $ad_q$, 
 $$ad_q = (\mu_{L} \otimes 
\mu_{R})(id \otimes S)\Delta,\eqno(2.7)$$  
where $\mu_{L}$ (respectively, $\mu_{R}$) 
is the left (respectively, right) multiplication: 
$\mu_{L}(x)y = xy$ (respectively, 
$\mu_{R}(x)y = (-1)^{degx.degy} yx$).
Then, the generators $E_{13}$ and $E_{31}$ defined 
in (2.3) can be written in an adjoint form,  
$$E_{13} = ad_q(E_{12})(E_{23}q^{H_2}). q^{H_1-H_2}~,\qquad 
E_{31} = - ad_q(E_{21})(E_{32}),\eqno(2.8)$$ 
being a $q$--analogue of  the classical one,  
$$e_{13}=ad(e_{12})e_{23}\equiv  [e_{12},e_{23}],\quad 
e_{31}=-ad(e_{21})e_{32}\equiv -[e_{21},e_{32}].$$ 
Of course, one can rescale the generators $E_{12}$ and $E_{23}$ to make $E_{13}$ 
in (2.8) to resemble more its classical counterpart $e_{13}=[e_{12},e_{23}]$.

    We see from the relations (2.1) -- (2.3) that each of the odd spaces $A_{\pm}$, 
$$A_+= \mbox{lin. env} .\{E_{13}, E_{23}\},\eqno (2.9)$$ 
$$A_-= \mbox{lin. env} .\{E_{31}, E_{32}\},\eqno (2.10)$$ 
is a representation space of the even subalgebra $U_q[gl(2/1)_0]\equiv 
U_q[gl(2)\oplus gl(1)]$, which, generated by generators $E_{12}$, 
$E_{21}$,  and $E_{ii}$, $i=1, 2, 3$, is a stability subalgebra of 
$U_q[gl(2/1)]$. Therefore,  we can construct a representation of 
$U_q[gl(2/1)]$ induced from some (finite--dimensional irreducible, 
for example) representation of $U_q[gl(2/1)_0]$ which is realized in 
a representation space (module) $V^{q}$ being a tensor product of a 
$U_q[gl(2)]$-module $V^q(gl_2)$ and a $gl(1)$-module ($gl(1)$-factors) 
$V^q(gl_1)$.  Let us take  throughout this paper $V^{q}$  to be an irreducible 
(later also finite--dimensional) $U_q[gl(2/1)_0]$-module. If we demand 
$$E_{23}V^{q}=0\eqno (2.11)$$
hence 
$$U_q(A_+)V^{q}=0,\eqno (2.12)$$
we turn the $U_q[gl(2/1)_0]$-module $V^{q}$ into a $U_q(B)$-module 
with  
$$B=A_+\oplus gl(2)\oplus gl(1).\eqno (2.13)$$
The $U_q[gl(2/1)]$-module $W^{q}$ induced from  the $U_q[gl(2/1)_0]$-
module $V^{q}$ is the factor space 
$$W^{q}= [U_q\otimes V^{q}]/I^{q}\eqno (2.14)$$ 
where 
$$U_q\equiv U_q[gl(2/1)],\eqno (2.15)$$ 
and $I^{q}$ is the subspace 
$$I^{q}=\mbox{lin. env} .\{ub\otimes v - u\otimes bv \| u\in U_{q }, b\in U_q(B)
\subset U_q, v\in V^{q}\}.\eqno (2.16)$$
By construction, any vector $w\in W^{q}$ can be represented as  
$$w= u\otimes v,~~~u\in U_q,~~~v\in V^{q}.\eqno (2.17)$$ 
Then $W^{q}$ is a $U_q[gl(2/1)]$-module in the sense 
$$gw\equiv g(u\otimes v)=gu\otimes v\in W^{q}\eqno (2.18)$$ 
for $g, u\in U_q, w \in W^{q}$ and $v \in V^{q}$.

  Using the commutation relations (2.1) -- (2.2) and the definitions (2.3) we
can prove an analogue of the Poincar\'e--Birkhoff--Witt theorem.\\[5mm]
{\bf Proposition 1}: {\it The quantum deformation $U_q:=U_q[gl(2/1)]$ 
is spanned on all possible linear combinations of the elements
$$ g= (E_{23})^{\eta_1}(E_{13})^{\eta_2}(E_{31})^{\theta_1}(E_{32})^
{\theta_2}g_0,\eqno (2.19)$$ 
where ~ $\eta_i, \theta_i=0, 1$ and 
$g_0 \in U_q[gl(2/1)_0]\equiv U_q[gl(2)]\oplus gl(1)].$}~
\\[7mm] 
Then the following proposition can be also proved :\\[5mm] 
{\bf Proposition 2}: {\it The induced $U_q[gl(2/1)]$-module $W^q$ is
the  linear span}
$$W^q=\mbox{lin. env.}\{(E_{31})^{\theta_1}(E_{32})^{\theta_2}\otimes v\| 
v \in V^{q}, ~~\theta _1, \theta_2 = 0, 1 \},\eqno(2.20)$$
{\it and, consequently, the set of all the vectors} 
$$\left|\theta_1,\theta_2; (m)\right>:= (E_{31})^{\theta_1}(E_{32})^{\theta_2}\otimes (m), 
~~ \theta_1, \theta_2=0,1,\eqno(2.21)$$ 
{\it constitutes a basis of $W^q$,
with $(m)$ being a basis of $V^q$.}\\[4mm] 
Thus, we can write the
$U_q[gl(2/1)]$-module $W^q$ in the form   
$$W^q([m])= T^q\otimes
V^q([m]),\eqno(2.22)$$ where 
$$T^q=\mbox{lin. env.}\left\{(E_{31})^{\theta_1}(E_{32})^{\theta_2},
\theta_i=0,1\right\}\eqno(2.23)$$ 
and $[m]$ is a \underline{signature} (an highest weight, in the case of
finite--dimensional representations) characterizing the  module $V^q$ 
and, therefore, also the module $W^q$. The basis (2.21) referred to as the 
\underline{induced basis} of $W^q$ is a tensor product  
$$\left|\theta_1,\theta_2; (m)\right>= \left|\theta_1,\theta_2\right>\otimes
(m)\eqno(2.21')$$ 
between a basis  
$$\left|\theta_1,\theta_2\right> := (E_{31})^{\theta_1}(E_{32})^{\theta_2},
~~ \theta_i =0,1,\eqno(2.24)$$ 
of $T^q$ and a basis  $(m)$ of $V^q$.

  Taking the fact that  $$ad_q\left (U_q[gl(2/1)_0]\right)T^q
\subset T^q\eqno(2.25)$$ 
we can consider $T^q$ as a module of the even subalgebra $U_q[gl(2/1)_0]$. 
This module is completely reducible since it represents a direct sum of three 
irreducible submodules  
$$T^q = T^q_0\oplus T^q_1\oplus T^q_2,\eqno(2.26)$$ where
\begin{eqnarray}
~~~~~~~~~~~~~~~~~~~~~~~~~~~~~~~
T^q_0&=&\mbox{lin. env.}\left\{(E_{31})^0(E_{32})^0\equiv {\bf 1}\right\}
\equiv {\bf C}, 
~~~~~~~~~~~~~~~~~ (2.27a)\nonumber \\ 
T^q_1&=&\mbox{lin. env.}\left\{E_{31},E_{32}\right\},
~~~~~~~~~~~~~~~~~~~~~~~~~~~~~~~~~~~ (2.27b)\nonumber \\ 
T^q_2&=&\mbox{lin. env.}\left\{E_{31}E_{32}\right\}.\nonumber 
~~~~~~~~~~~~~~~~~~~~~~~~~~~~~~~~~~~~~ (2.27c)
\end{eqnarray} 
 Every subspace $T^q_i$, $i=0,1,2$, as an irreducible $U_q[gl(2/1)_0]$-module, 
is characterized by a signature, say $[\mu]_i$, which is always fixed and will be 
determined in the next section (see (3.11) ): 
$$T^q_i=T^q_i([\mu]_i).\eqno(2.28)$$ 
So, the module $W^q$ being a tensor product of two $U_q[gl(2/1)_0]$-modules, 
$T^q$ and  $V_0^q$, is also a $U_q[gl(2/1)_0]$-module which, in general, is 
reducible and can be written now in the form 
$$W^q([m])= D_0\oplus D_1\oplus D_2,\eqno(2.29)$$ 
where 
$$D_i= T^q_i\otimes V^q([m]), ~~ i=0,1,2.\eqno(2.30)$$
Here, as seen later, $D_0$ and $D_2$ are irreducible $U_q[gl(2/1)_0]$-modules, 
but $D_1$ is a reducible one (see, (3.23) ).\\[5mm] 
{\bf Proposition 3}: 
{\it The $U_q[gl(2/1)]$-module 
$W^{q}$ is decomposed into (four or less) finite--dimentional irreducible modules 
$V_k^{q}$ of the even subalgebra $U_q[gl(2/1)_0]$, 
$$W^{q}([m])= \bigoplus_{0\leq k\leq 3}V_k^{q}([m]_k),\eqno(2.31)$$ 
where $[m]$ and $[m]_k$ are some signatures (highest-weights) characterizing 
the module $W^{q}\equiv W^{q}([m])$ and the modules 
$V^{q}_k\equiv V_k^{q}([m_k])$, respectively.}

  Now we are ready to construct finite--dimensional representations of
$U_q[gl(2/1)]$ in a basis of its even subalgebra $U_q[gl(2/1)_0]$. 
These representations are induced from finite--dimensional irreducible
representations of the even subalgebra $U_q[gl(2/1)_0]$. For a basis 
of the latter we can chose a Gel'fand--Zetlin (GZ) one.\\[1cm]
{\bf III. FINITE -- DIMENSIONAL REPRESENTATIONS OF $U_q[gl(2/1)]$}
 
  A finite--dimensional representation of $U_q[gl(2/1)_0]$ is realized 
in some space (module) which could be one of the above $V_k^{q}$ 
whose basis, a $U_q[gl(2/1)_0]$-basis,  can be chosen as a tensor product 
$$\left[ 
 \begin{array}{c}
m_{12}~~~m_{22}\\
 m_{11}
 \end{array};
\begin{array}{c}
m_{32}
\\
m_{31}
\end{array}\right]\equiv (m)_{gl(2)}\otimes m_{31}\equiv (m)_k
\eqno (3.1a)$$\\
between a (GZ) basis  $(m)_{gl(2)}$ of $U_q[gl(2)]$ and $gl(1)$-factors  
$m_{31}$, where $m_{ij}$ are complex numbers such that
$$m_{12}-m_{11}, ~ m_{11}-m_{22} \in {\bf Z}_+,\eqno (3.1b)$$
$$m_{32}=m_{31}.\eqno (3.1c)$$
Indeed, finite--dimensional representations of  $U_q[gl(2)]$ are  
highest weight and the generators $E_{ij}$, $i,j=1,2$, and $E_{33}$ 
(called the even generators of $U_q[gl(2/1)]$) really satisfy the 
commutation  relations  (2.1a) -- (2.1d)  for $U_q[gl(2/1)_0]$ if 
they are defined on (3.1) as follows 
\begin{eqnarray*}
  ~~~~~~~~~~~~~~~~~~~~E_{11}(m)_k
&=&(l_{11}+1)(m)_k,\\
 E_{22}(m)_k &=& (l_{12}+l_{22}-l_{11}+2)(m)_k,\\
E_{12}(m)_k &=&([l_{12}-l_{11}][l_{11}-l_{22}])^{1/2}(m)_k^{+11},\\
E_{21}(m)_k &=& ([l_{12}-l_{11}+1][l_{11}-l_{22}-1])^{1/2}(m)_k^{-11},\\
E_{33}(m)_k
&=&(l_{31}+1)(m)_k,~~~
~~~~~~~~~~~~~~~~~~~~~~~~~~~~~~~~~~~~~~~~~~~~~~~\label{}{(3.2)}
\end{eqnarray*}
where $$l_{ij}=m_{ij}-(i-2\delta_{i,3}),\eqno (3.3)$$
and $(m)_k^{\pm ij}$ is a vector obtained from $(m)$ by replacing $m_{ij}$ 
with $m_{ij}\pm 1$. The signature of a basis vector $(m)_k$ now is the highest 
weight described by the first (top) row  of the patterns (3.1)
$$[m]_k= [m_{12}, m_{22}, m_{32}]\eqno (3.4)$$
remaining unchanged under the action of the even generators is nothing but an
ordered set of eigen-values of the Cartan
generators $E_{ii}$, $i=1,2,3$, on
the highest weight vector $(M)_k$ defined
as follows
$$E_{12}(M)_k=0,\eqno(3.5)$$
$$E_{ii}(M)_k=m_{i2}(M)_k.\eqno(3.6)$$ 
The highest weight vector $(M)_k$ is a vector $(m)_k$  with $m_{11}$ taking the 
maximal value $m_{11}=m_{12}$, 
$$(M)_k=\left[ 
\begin{array}{c}
m_{12}~~~m_{22}\\
m_{12}
\end{array} ;
\begin{array}{c}
m_{32}=m_{31}\\
m_{31}
\end{array}\right],\eqno (3.7)$$\\
and vice versa a (lower weight) vector $(m)_{k}$ can be derived from $(M)_{k}$ 
via the formula 
\begin{eqnarray*}
~~~~~~~~~~~~~~~~(m)_{k}&=&\left({[m_{11}-m_{22}]!\over
[m_{12}-m_{22}]![m_{12}-m_{11}]!}\right)^{1/2}
(E_{21})^{m_{12}-m_{11}}(M)_{k}.
~~~~~~~~~~~~(3.8) 
\end{eqnarray*} 
\\
The subscript $k$ in the l.h.s of (3.4) can be omitted when there is no degeneration 
among signatures of basis vectors. Additionally, for the case $k=0$, as will be 
seen $V_0^{q} \equiv V^{q}$, we can always skip the subscript 0, 
$$(m)_0\equiv (m),~~~[m]_0\equiv [m],~~~(M)_0\equiv (M).\eqno (3.9)$$

   In a GZ basis (3.1),  the highest weights $[\mu]_i$ of the subspaces $T^q_i$  
have the form  (3.4), $[\mu]_i\equiv [\mu_{12},\mu_{22},\mu_{32}]$, that is  
$$T^q_i=T^q_i([\mu]_i)\equiv T^q_i([\mu_{12},\mu_{22},\mu_{32}]).$$ 
Let us denote the GZ basis vectors of $T^q_i([\mu]_i)$ by 
$$\left[   \begin{array}{c} 
\mu_{12}~~~\mu_{22}\\  \mu_{11}  \end{array}; 
\begin{array}{c} 
\mu_{32}=\mu_{31}\\  
\mu_{31} 
\end{array}\right]\equiv (\mu)_{gl(2)}\otimes \mu_{31}\equiv (\mu).\eqno (3.10)$$
\\ 
  Using the action of $U_q[gl(2/1)_0]$ on $T^q_i$ we identify 
the basis vectors (2.24) as follows:  
$$\left|0,0\right>\equiv 1= \left[ 
\begin{array}{c} 
0~~~0\\  
0 
\end{array}; 
\begin{array}{c} 
0\\ 
0 
\end{array}\right]~ \in~ T^q([0,0;0])=T^q_0,\eqno (3.11a)$$ 

$$
\left.
\begin{array}{lll}
\left|1,0\right>&\equiv &E_{31}= -\left[ 
\begin{array}{c} 
0~~~-1\\ 
-1\end{array}; 
\begin{array}{c} 
1\\
1\end{array}\right]  
\\[9mm] 
\left|0,1\right>&\equiv &E_{32}= \left[ 
\begin{array}{c} 
0~~~-1\\ 
0 
\end{array}; 
\begin{array}{c} 
1\\1 
\end{array}\right] 
\end{array}
\right\} ~ \in~ T^q([0,-1;1])=T^q_1,
\eqno(3.11b)
$$

$$\left|1,1\right>\equiv  E_{31}E_{32}=\left[ 
\begin{array}{c} 
-1~~~-1\\ 
-1 
\end{array}; 
\begin{array}{c} 
2\\2 
\end{array}\right]~ \in~ T^q([-1,-1;2])=T^q_2. 
\eqno (3.11c)$$ 
\\
In the latter formulae the subscripts $i$ of signatures of the subspaces $T^q_i$ 
can be omitted since there is no degeneration among these signatures. 

 We can combine all the basis vectors (3.11) in a common formula: 
$$\left | \theta_1, \theta_2 \right >
= (-1)^{\theta_1 (\theta_2 +1)}(\mu), \eqno(3.12)$$ 
where $\theta_1, \theta_2 =0,1$,  and 
\begin{eqnarray*}
~~~~~~~~~~~~~~~~~~~~~~ (\mu)&=&
\left[  
\begin{array}{c} 
\mu_{12} ~~~~~ \mu_{22}\\[5mm] \mu_{11}  
\end{array}; 
\begin{array}{c} 
\mu_{32}\\[5mm] \mu_{31}
\end{array}\right]\\[5mm]
&=&
\left[  \begin{array}{c} 
-\theta_1\theta_2 ~~~ -{1\over 2}\left\{1 +(-1)^{(1-\theta_1)
(1-\theta_2)}\right\}\\[5mm] -\theta_1  
\end{array}; 
\begin{array}{c} 
\theta_1 + \theta_2 \\[5mm] \theta_1 + \theta_2
\end{array}\right].
~~~~~~ ~ (3.13)
\end{eqnarray*} \\
The action of the even subalgebra $U_q[gl(2/1)_0]$ 
on the basis (3.11) of  $T^q_i$  is the following  
\begin{eqnarray*}
~~~~~~~~~~~~~~~~~~  E_{ij}\left | \theta_1, \theta_2 \right > &=& 
-\theta_i(1-\theta_j).\left |1-\theta_1, 1-\theta_2 \right >,~~ \mbox{for} 
~~ i,j=1,2, ~~ i\neq j,\\
E_{ii}\left | \theta_1, \theta_2 \right > &=& 
- \theta_i.\left |\theta_1, \theta_2 \right >, ~~ i=1,2,\\
E_{33}\left | \theta_1, \theta_2 \right > &=& 
(\theta_1+\theta_2).\left |\theta_1, \theta_2 \right >. ~~~~~~~~~~~~~~ 
~~~~~~~~~~~~~~~~~~~~~~~~~~~~~ (3.14)   
\end{eqnarray*}
Now the induced basis (2.21) can be written in the form 
$$ \left|\theta_1,\theta_2; (m)\right>= (-1)^{\theta_1 (\theta_2 +1)}(\mu)\otimes (m).
\eqno(3.15)$$
To find the transformation of the latter basis under $U_q[gl(2/1)]$, it is sufficient to 
find transformations of this basis under the Weyl--Chevalley generators, which are 
those $E_{ij}$ with $|i-j|\leq 1$, $i,j=1,2,3$.  The actions of the even generators follow 
from their co-product structure and their actions (3.2) on $(\mu)$ and $(m)$, while those 
of the odd generators follow from 
\begin{eqnarray*} 
~~~~~ E_{32}(E_{31})^{\theta_1}(E_{32})^{\theta_2}&=& 
(-q)^{\theta_1}(E_{31})^{\theta_1}(E_{32})^{\theta_2+1}, \\
E_{23}(E_{31})^{\theta_1}(E_{32})^{\theta_2}&=& 
(-1)^{\theta_1+\theta_2}(E_{31})^{\theta_1}(E_{32})^{\theta_2}E_{23}
+(-1)^{\theta_1}\theta_2(E_{31})^{\theta_1}[H_2]\\
&&-\theta_1\theta_2E_{31}q^{-H_2-1}+\theta_1q^{-\theta_2}(E_{32})^{\theta_2}
E_{21}q^{-H_2-1}.~~~~~~~~~~~~~~~ (3.16)
\end{eqnarray*} 
The latter in turn can follow from a more general (deformed) 
commutation relation
\begin{eqnarray*} 
E_{ij}(E_{31})^{\theta_1}(E_{32})^{\theta_2}&=&
q^{
(\delta_{i3}\delta_{j2}+\delta_{i2}\delta_{j1})
\theta_1-
\delta_{i2}\delta_{j1}
\theta_2}
(-1)^{
\delta_{i2}\delta_{j3}
(\theta_1 +\theta_2) +
\delta_{i3}\delta_{j2}
\theta_1}
(E_{31})^{\theta_1}(E_{32})^{\theta_2}E_{ij}\\[2mm]
&&-\left\{
\delta_{i1}\delta_{j1}
\theta_1+(
\delta_{i2}\delta_{j2}-
\delta_{i3}\delta_{j3})\theta_2\right\}(E_{31})^{\theta_1}(E_{32})^{\theta_2}\\[3mm]
&&-\theta_i(1-\theta_j)(E_{31})^{1-\theta_1}(E_{32})^{1-\theta_2}
q^{
\delta_{i1}\delta_{j2}
H_1}\\
&&+
\delta_{i2}\delta_{j3}
\theta_1\left\{q^{-\theta_2}
(E_{32})^{\theta_2}E_{21}-\theta_2
E_{31}\right\}q^{-1-H_2}\\
&&+
\delta_{i2}\delta_{j3}
\theta_2
(-E_{31})^{\theta_1}
[H_2], ~~~~~~~~~~~~~~~~~~~~~~~~~~~~~~~~~
~~~~~~~~~~~~~ (3.17) 
\end{eqnarray*}
where 
$i, j=1, 2, 3$, $|i-j|\leq 1$,  
 $ \theta_1, \theta_2=0, 1$, $\theta_3=\theta_2$.  
This  commutation relation is, of course, consistent with (3.14). 

  Taking into account (3.11) -- (3.17) we get representations of $U_q[gl(2/1)]$ 
in the induced basis (2.21)
\begin{eqnarray*} 
E_{ij}\left\|\theta_1,\theta_2;(m)\right >&=&(1-\delta_{i2}\delta_{j3})
q^{
(\delta_{i3}\delta_{j2}+\delta_{i2}\delta_{j1})
\theta_1-
\delta_{i2}\delta_{j1}
\theta_2}
(-1)^{
\delta_{i2}\delta_{j3}
(\theta_1 +\theta_2) +
\delta_{3i}\delta_{j2}
\theta_1}\\
&&\times
\left|\theta_1,\theta_2 +\delta_{i3}\delta_{j2}; 
(m)_{ij}\right >\\[2mm]
&&-\left\{
\delta_{i1}\delta_{j1}
\theta_1+(
\delta_{i2}\delta_{j2}-\delta_{i3}\delta_{j3}
)\theta_2\right\}
.\left|\theta_1,\theta_2;(m)\right >\\
&&-\theta_i(1-\theta_j)
q^{
\delta_{i1}\delta_{j2}
h_1}.\left|1-\theta_1,1-\theta_2;(m)\right >\\
&&+
\delta_{i2}\delta_{j3}\theta_1 q^{-1-h_2}\left\{
q^{-\theta_2}
.\left|1-\theta_1, \theta_2; (m)_{21}\right>-\theta_2. \left|
\theta_1, 1-\theta_2; (m)\right>\right\}\\
&&+
\delta_{i2}\delta_{j3}
\theta_2(-1)^{\theta_1}
[h_2]. \left|\theta_1, 1-\theta_2; (m)\right>, ~~~~~~~~~~~~~~~~~ 
~~~~~~~~~~~ (3.18a)
\end{eqnarray*}
where  $l$ and $h_i$ are 
respectively eigenvalues of $L$ and $H_i$ on $(m)$, while   
\begin{eqnarray*}
(m)_{ij}
&=&\left\{
\begin{array}{l}
 E_{ij}(m), \mbox{given by  (3.2)},  ~\mbox{if} ~~  i,j=1,2 ~~ \mbox{or} ~ ~ i=j=3,\\
(m) ~,  ~\mbox{otherwise} .
\end{array}
\right.  ~~~~~~~~~~~~~~~~~~ (3.18b)
\end{eqnarray*}
These transformations give different representations for different $[m]$. The 
representations of  $U_q[gl(2/1)]$ constructed are in general reducible. 
However, the induced basis is not convenient for investigating the representation 
structure. Let us go to another, more appropriate to this goal, basis.

  The module $V^q$ is a tensor product of a $U_q[gl(2)]$-module with a 
$U_q[gl(1)]$-module (in fact, a $gl(1)$-factor), 
$$V^q([m_{12},m_{22},m_{32}])=V^q([m_{12},m_{22}])\otimes 
V^q([m_{32}]),\eqno (3.19a)$$  
and so is the module $T^q$,   
$$T^q([\mu_{12},\mu_{22},\mu_{32}])=T^q([\mu_{12},\mu_{22}])\otimes 
T^q([\mu_{32}]).\eqno (3.19b)$$  
Then the module $W^q$ in (2.22) can be written as follows 
$$W^q([m]) = \left\{T^q([\mu_{12},\mu_{22}])\odot V^q([m_{12},m_{22}])\right\}
\otimes \left\{T^q([\mu_{32}])\odot V^q([m_{32}])\right\}.\eqno(3.20)$$
Here the notation $\odot $ is used for a tensor product between two modules of one 
and the same (quantum) algebra, whereas $\otimes$ is a more general notation used 
for a tensor product of two arbitrary spaces or modules. In general, the  
$U_q[gl(2)]$-module $T^q([\mu_{12},\mu_{22}])\odot V([m_{12},m_{22}])$ in 
(3.20) is reducible and can be decomposed into a direct sum of irreducible modules 
$$T^q([\mu_{12},\mu_{22}])\odot V^q([m_{12},m_{22}])=\bigoplus_{i=0}^{n} 
V^q([\mu_{12}+m_{12}-i,\mu_{22}+m_{22}+i]),\eqno (3.21)$$ 
where $$n= \mbox{min}~(\mu_{12}-\mu_{22},m_{12}-m_{22}),$$ 
\\
while the $gl(1)$-factor $T^q([\mu_{32}])\odot V^q([m_{32}])$ is just 
$$T^q([\mu_{32}])\odot V^q([m_{32}])=V^q([\mu_{32}+m_{32}]).\eqno (3.22)$$ 
Taking into account (2.30) and (3.19)  -- (3.22) we get 
\begin{eqnarray}
~~~~~~~~~~~~~~~~~~~~~~~~~
D_0&\equiv & T^q([0,0,0])\odot V^q([m_{12},m_{22},m_{32}])
\equiv  V^q([m_{12},m_{22},m_{32}]),\nonumber \\[3mm]
D_1&\equiv &T^q([0,-1,1])\odot V^q([m_{12},m_{22},m_{32}])\nonumber \\
&= &\bigoplus_{i=0}^1 V^q([m_{12}-i,m_{22}+i-1,m_{32}+1]),\nonumber \\[3mm]
D_2&\equiv & T^q([-1,-1,2])\odot V^q([m_{12},m_{22},m_{32}])\nonumber \\
&= &V^q([m_{12}-1,m_{22}-1,m_{32}+2]).
~~~~~~~~~~~~~~~~~~~~~~~~~~~~ (3.23)
\nonumber  
\end{eqnarray} 
Inserting (3.23) in (2.28) we prove (2.31) with $V^q_k$ identified as 
follows 
\begin{eqnarray} 
~~~~~~~~~~~~~~~~~~~~~~~~~ 
V^q_0&\equiv & V^q([m_{12},m_{22},m_{32}]) =V^q,\nonumber \\[3mm]
V^q_1&\equiv &V^q([m_{12},m_{22}-1,m_{32}+1]),\nonumber \\[3mm]
V^q_2&\equiv &V^q([m_{12}-1,m_{22},m_{32}+1]),\nonumber \\[3mm]
V^q_3&\equiv &V^q([m_{12}-1,m_{22}-1,m_{32}+2]).
~~~~~~~~~~~~~~~~~~~~~~~~~~~~ (3.24)
\nonumber  
\end{eqnarray}

  Instead of the induced basis (2.21) for a basis of $W^q$ we can chose the 
union of the bases, the GZ bases (3.1) in the case, of all its subspaces $V^q_k$. 
This  new basis of  $W^q$ is referred to as its \underline{reduced basis} which 
is related to the induced one (2.21) via the Clebsch--Gordan (CG) decomposition. 
In order to derive  such a relation between the two bases for the whole $W^q$ we 
should have it first for each of the subspaces (3.21) and (3.22).  Within the subspace 
(3.21), which is a $U_q[gl(2)]$-module,  the relation between 
the induced basis 
$$(\mu)_{gl(2)}\odot (m)_{gl(2)}\equiv \left[\begin {array}{c} 
\mu_{12}~~\mu_{22}\\ 
\mu_{11}\end{array}\right]\odot \left[ 
\begin{array}{c} 
m_{12}~~m_{22}\\ 
m_{11}\end{array}\right]~ \in ~ T^q([\mu_{12},\mu_{22}])\odot  
V^q([m_{12},m_{22}]) $$\\
and the reduced basis  
$$(m')_{gl(2)}\equiv \left[ \begin{array}{c} 
m'_{12}~~~m'_{22}\\ 
m'_{11}\end{array}\right]~ \in ~ V^q([m'_{12},m'_{22}]),$$ 
$$m'_{12}=\mu_{12}+m_{12}-i,~ m'_{22}=\mu_{22}+m_{22}+i,$$  
can be written in the form 
$$(m')_{gl(2)}=\sum_{\mu_{11},m_{11}} 
\left[
\begin{array}{c} 
m'_{12}~~m'_{22}\\ 
m'_{11} 
\end{array}
\right| 
\left. 
\begin{array}{c}\mu_{12}~~\mu_{22}\\ 
\mu_{11}\end{array}; 
\begin{array}{c} 
m_{12}~~m_{22}\\ 
m_{11}\end{array}
\right] 
(\mu)_{gl(2)}\odot (m)_{gl(2)}
~~, \eqno (3.25)$$ 
\\where
$$\left[
\begin{array}{c} 
m'_{12}~~m'_{22}\\ 
m'_{11} 
\end{array}
\right| 
\left. 
\begin{array}{c}\mu_{12}~~\mu_{22}\\ 
\mu_{11}\end{array}; 
\begin{array}{c} 
m_{12}~~m_{22}\\ 
m_{11}\end{array}
\right] 
\eqno (3.26)$$ 
\\are the  Clebsch--Gordan coefficients of $U_q[gl(2)]$. 
The relation between the two bases within the subspace (3.22) is simply 
$$m'_{31}=\mu_{31}+m_{31}.\eqno(3.27)$$  
Now taking into account (3.25) -- (3.27) we can express the reduced 
basis of $W^q$ in terms of the induced one (2.21) :

\begin{eqnarray*} 
~~ (m)_0&\equiv & \left[ 
\begin{array}{c} 
m_{12} ~~~ m_{22}\\ 
m_{11} 
\end{array}; 
\begin{array}{c} 
m_{32}\\ 
m_{32} 
\end{array} 
\right]\equiv (m)  ~ \in ~ V_0^q\equiv
V^q,\\[5mm] 
(m)_1&\equiv & \left[ 
\begin{array}{c} 
m_{12} ~~~ m_{22}-1\\ 
m_{11} 
\end{array}; 
\begin{array}{c} 
m_{32}+1\\ 
m_{32}+1 
\end{array} 
\right]\\[3mm] 
&=& \left[ 
\begin{array}{c} 
m_{12} ~~ m_{22}-1\\ 
m_{11} 
\end{array} 
\right| 
\left. 
\begin{array}{c} 
0~~-1\\ 
0  
\end{array}; 
\begin{array}{c} 
m_{12} ~~ m_{22}\\ 
m_{11} 
\end{array} 
\right].\left|0,1; (m)\right>\\[3mm] 
& &-\left[ 
\begin{array}{c} 
m_{12} ~~ m_{22}-1\\ 
m_{11} 
\end{array}
\right| 
\left. 
\begin{array}{c} 
0~~-1\\ 
-1 
\end{array}; 
\begin{array}{c} 
m_{12} ~~ m_{22}\\ 
m_{11}+1 
\end{array} 
\right].\left|1,0; (m)^{+11}\right> ~ \in ~ V_1^q,\\[5mm] 
(m)_2&\equiv &  
\left[ 
\begin{array}{c} 
m_{12}-1 ~~~ m_{22}\\ 
m_{11} 
\end{array}; 
\begin{array}{c} 
m_{32}+1\\ 
m_{32}+1 
\end{array} 
\right]\\[3mm] 
&=&\left[ 
\begin{array}{c} 
m_{12}-1 ~~ m_{22}\\ 
m_{11} 
\end{array} 
\right| 
\left. 
\begin{array}{c} 
0~~-1\\ 
0 
\end{array}; 
\begin{array}{c} 
m_{12} ~~ m_{22}\\ 
m_{11} 
\end{array} 
\right]. \left|0,1; (m)\right>\\[3mm] 
& &-\left[ 
\begin{array}{c} 
m_{12}-1 ~~ m_{22}\\ 
m_{11} 
\end{array} 
\right| 
\left. 
\begin{array}{c} 
0~~-1\\ 
-1 
\end{array}; 
\begin{array}{c} 
m_{12} ~~ m_{22}\\ 
m_{11}+1 
\end{array} 
\right].\left|1,0; (m)^{+11}\right> ~ \in ~ V_2^q,\\[5mm] 
(m)_3&\equiv &\left[ 
\begin{array}{c} 
m_{12}-1 ~~~~ m_{22}-1\\ 
m_{11} 
\end{array};  
\begin{array}{c} 
m_{32}+2\\ 
m_{32}+2 
\end{array} 
\right]\\[2mm] 
&=&\left[ 
\begin{array}{c} 
m_{12}-1 ~~ m_{22}-1\\ 
m_{11} 
\end{array} 
\right| \left. 
\begin{array}{c} 
-1~~-1\\ 
-1 
\end{array}; 
\begin{array}{c} 
m_{12} ~~ m_{22}\\ 
m_{11}+1 
\end{array} 
\right].\left|1,1; (m)^{+11}\right> ~ \in ~ V_4^q.~~~~~\label{}{(3.28)}  
\end{eqnarray*} 
\\
 By construction, $W^q$ is characterized by the  highest weight  of  $V^q$, 
the signature $[m]$  in (3.9). 
In order to describe $W^q$ as the whole we unify the basis vectors (3.28) in  
a single notation,  
$$
\left[ 
\begin{array}{c}
m_{13}~~~m_{23}~~~m_{33}\\ 
m_{12}~~~m_{22}~~~m_{32}\\ 
m_{11}~~~~0~~~m_{31} 
\end{array}\right],
\eqno (3.29)$$ 
\\ 
by putting on the top of the GZ patterns (3.28) an additional row 
which is exactly the highest weight $[m]$ of $V^q$, denoted now as 
$$[m]\equiv [m_{13},m_{23},m_{33}].\eqno(3.30)$$
This row of (3.29) remains unchanged throughout the whole $W^q$, while the 
second row depending on $k$ represents the first row of one of the patterns (3.28) 
and tells us which subspace $V^q_k$ the considered basis vector (3.29)  of  $W^q$ 
belongs to.  The basis (3.29) reflects the branching rule 
$U_q[gl(2/1)]\supset U_q[gl(2)\otimes gl(1)]$ and it can be called a 
quasi-GZ basis. The subspaces (3.24) in this new notation is 
\begin{eqnarray} 
~~~~~~~~~~~~~~~~~~~~~~~~~~~~~
V^q_0&\equiv &V^q([m_{13},m_{23},m_{33}]) \equiv V^q,\nonumber \\[3mm]
V^q_1&\equiv &V^q([m_{13},m_{23}-1,m_{33}+1]),\nonumber \\[3mm]
V^q_2&\equiv &V^q([m_{13}-1,m_{23},m_{33}+1]),\nonumber \\[3mm]
V^q_3&\equiv &V^q([m_{13}-1,m_{23}-1,m_{33}+2]).
~~~~~~~~~~~~~~~~~~~~~~~ (3.31) 
\nonumber  
\end{eqnarray} 

  Let us now determine the CG coefficients in (3.28).  To do that we use the Hopf 
algebra structure which is again helpful. We start with the subspace $V^q_1$. The 
highest vector here is 
$$(M)_1=a_1 \left(E_{32}\otimes (M)\right),\eqno (3.32)$$ 
where $a_1$ is an arbitrary complex coefficient which may depend on $q$. 
Formula (3.8) now becomes  
$$(m)_1=\left(\frac {[l_{11}-l_{23}]!}{[2l+1]![l_{13}-l_{11}]!}\right)^{1/2} 
(E_{21})^{l_{13}-l_{11}}(M)_1, \eqno (3.33a)$$ 
where 
$$l=(m_{13}-m_{23})/2.\eqno (3.33b)$$ 
Replacing 
\begin{eqnarray*} 
~~~~~~~~
(E_{21})^{l_{13}-l_{11}}(M)_1&\equiv&  a_1\left\{\Delta
(E_{21})\right\}^{l_{13}-l_{11}}\left(E_{32}\otimes (M)\right)\\ 
&=&-a_1[l_{13}-l_{11}]\left( 
\frac{[2l]![l_{13}-l_{11}-1]!}{[l_{11}-l_{23}]!}\right)^{1/2}\left(E_{31}
\otimes (m)^{+11}\right)\\  
&&+ a_1q^{l_{11}-l_{13}}\left( 
\frac{[2l]![l_{13}-l_{11}]!}{[l_{11}-l_{23}-1]!}\right)^{1/2}\left(E_{32}
\otimes (m)\right) 
~~~~~~~~~~~~ \label{}{(3.34)} 
\end{eqnarray*} 
in (3.33) we obtain  
$$(m)_1=a_1\left\{-\left( 
\frac{[l_{13}-l_{11}]}{[2l+1]}\right)^{1/2}.\left|1,0; (m)^{+11}\right>
+q^{l_{11}-l_{13}}\left( 
\frac{[l_{11}-l_{23}]}{[2l+1]}\right)^{1/2}.\left|0,1; (m)~\right>\right\}. $$\\
So with the help of the Hopf algebra structure the necessary CG coefficients can 
be found directly and easily without knowing in advance a general formula for  them.  
In the same way, from the highest weight vectors 
$(M )_2\in V^q_2$ and 
$(M )_3\in V^q_3$,  
$$(M)_2=a_2\{E_{31}\otimes (M)+q^{2l}[2l]^{-1/2}E_{32}\otimes (M)^{-11}\},$$ 
$$(M)_3=a_3E_{31}E_{31}\otimes (M),$$ 
we can find explicit expressions for $(m)_2$ and $(m)_3$, respectively. Thus,  we have 
the following relation between the reduced and the induced basis 
\begin{eqnarray*} 
~~~~~~ (m)_0&\equiv & \left[ 
\begin{array}{c} 
m_{13}~~~m_{23}~~~m_{33}\\ 
m_{13}~~~m_{23}~~~m_{33}\\ 
m_{11}~~~~0~~~~~m_{33} 
\end{array}\right]=
\left|0,0; (m)\right>\equiv (m),\\[4mm]
(m)_1&\equiv& \left[ 
\begin{array}{l} 
m_{13}~~~~~m_{23}~~~~~~~~m_{33}\\ 
m_{13}~~~m_{23}-1~~~m_{33}+1\\ 
m_{11}~~~~~~0~~~~~~~~~m_{33}+1 
\end{array}\right]\\[2mm] 
&=&a_1\left\{-\left( 
\frac{[l_{13}-l_{11}]}{[2l+1]}\right)^{1/2}.\left|1,0; (m)^{+11}\right>
+q^{l_{11}-l_{13}}\left( 
\frac{[l_{11}-l_{23}]}{[2l+1]}\right)^{1/2}.\left|0,1; (m)~\right>\right\},\\[4mm]
(m)_2&\equiv& \left[ 
\begin{array}{l} 
~~m_{13}~~~~~m_{23}~~~~~~m_{33}\\ 
m_{13}-1~~m_{23}~~~m_{33}+1\\ 
~~m_{11}~~~~~~0~~~~~~m_{33}+1 
\end{array}\right]\\[2mm] 
&=&a_2\left\{ 
\left( 
\frac{[l_{11}-l_{23}]}{[2l]}\right)^{1/2}.\left|1,0; (m)^{+11}\right>
\right.\\ &&\left.
+q^{l_{11}-l_{23}} 
\left(\frac{[l_{13}-l_{11}]}{[2l]}\right)^{1/2}.\left|0,1; (m)~\right>\right\},\\[4mm]
(m)_3&\equiv& \left[ 
\begin{array}{l} 
~~m_{13}~~~~~~~~m_{23}~~~~~~~m_{33}\\ 
m_{13}-1~~m_{23}-1~~m_{33}+2\\ 
~~m_{11}~~~~~~~~0~~~~~~~~m_{33}+2 
\end{array}\right]=a_3.\left|1,1; (m)^{+11}\right>  
~~~~~~~~~~~~~~~~~ (3.35)
\end{eqnarray*} 
\\
and, equivalently, the inverse relation 
\begin{eqnarray*} 
\left|0,0; (m)~\right>&=& (m)\\[4mm]
~~~~~~~~~~~~~~\left|1,0; (m)~\right>&=&-\frac{1}{a_1}q^{l_{11}-l_{23}-1} 
\left( 
\frac{[l_{13}-l_{11}+1]}{[2l+1]}\right)^{1/2}(m)_1^{-11}\\[4mm] 
&&+\frac{1}{a_2}q^{l_{11}-l_{13}-1} 
\frac{([l_{11}-l_{23}-1][2l])^{1/2}}{[2l+1]}(m)_2^{-11},\\[4mm] 
\left|0,1; (m)~\right>&=&\frac{1}{a_1} 
\left( 
\frac{[l_{11}-l_{23}]}{[2l+1]}\right)^{1/2}(m)_1\\[4mm] 
&&+\frac{1}{a_2} 
\frac{([l_{13}-l_{11}][2l])^{1/2}}{[2l+1]}(m)_2,\\[4mm] 
\left|1,1; (m)~\right>&=&\frac{1}{ a_3} (m)_3^{-11}.
~~~~~~~~~~~~~~~~~~~~~~~~~~~~~~~~~~~~~~~~~~~~~~~~~~~~~~~~  
\label{}{(3.36)} 
\end{eqnarray*} 

   Now we are ready to compute all the matrix elements of the generators in 
the basis (3.35) which allows a clear description of the structure of the  
module $W^{q}$.  Since the finite--dimensional representations of the  
$U_q[gl(2/1)]$ in some basis are completely defined by the actions of 
the even generators and the odd Weyl--Chevalley ones $E_{23}$ and $E_{32}$ 
in the same  basis, it is sufficient to write down the matrix elements of these 
generators only. For the even generators the matrix elements have already been 
given in  (3.2), while for $E_{23}$ and $E_{32}$, using the relations  (2.1)--(2.3), 
(3.35) and (3.36) we get 
\begin{eqnarray*} 
~~~~~~~~~~~~~~
E_{23}(m)&=&0,\\[4mm] 
E_{23}(m)_1&=&a_1 
\left( 
\frac{[l_{11}-l_{23}]}{[2l+1]}\right)^{1/2}[l_{23}+l_{33}+3](m),\\[4mm] 
E_{23}(m)_2&=&a_2 
\left( 
\frac{[l_{13}-l_{11}]}{[2l]}\right)^{1/2}[l_{13}+l_{33}+3](m),\\[4mm] 
E_{23}(m)_3&=&a_3 
\left\{\frac{1}{a_1q} 
\left( 
\frac{[l_{13}-l_{11}]}{[2l+1]}\right)^{1/2}[l_{13}+l_{33}+3](m)_1\right.\\[4mm] 
&&-\left.\frac{1}{a_2q}([l_{11}-l_{23}][2l])^{1/2} 
\frac{[l_{23}+l_{33}+3]}{[2l+1]}(m)_2\right\},\\[6mm] 
E_{32}(m)&=&\frac{1}{a_1} 
\left( 
\frac{[l_{11}-l_{23}]}{[2l+1]}\right)^{1/2}(m)_1 
+\frac{1}{a_2} 
\frac{([l_{13}-l_{11}][2l])^{1/2}}{[2l+1]}(m)_2,\\[4mm] 
E_{32}(m)_1&=&\frac{a_1q}{a_3} 
\left(\frac{[l_{13}-l_{11}]}{[2l+1]}\right)^{1/2}(m)_3,\\[4mm] 
E_{32}(m)_2&=&-\frac{a_2q}{a_3} 
\left( 
\frac{[l_{11}-l_{23}]}{[2l]}\right)^{1/2}(m)_3,\\[4mm] 
E_{32}(m)_3&=&0. 
~~~~~~~~~~~~~~~~~~~~~~~~~~~~~~~~~~~~~~~~~~~~
~~~~~~~~~~~~~~~~~~~~~~~~~~~  (3.37) 
\end{eqnarray*} 
 All the matrix elements of the Chevalley generators obtained here 
coincide, of course, with the ones obtained previously by another (but 
longer) way \cite{pq96, ic96}.  Besides that, we can easily find matrix 
elements for non--Chevalley generators too:
\begin{eqnarray*} 
~~~~~~~~~~~~  E_{31}(m)&=&-\frac{1}{a_1}q^{l_{11}-l_{23}-1} 
\left( 
\frac{[l_{13}-l_{11}+1]}{[2l+1]}\right)^{1/2}(m)_1^{-11}\\[4mm] 
&&+\frac{1}{a_2}q^{l_{11}-l_{13}-1} 
\frac{([l_{11}-l_{23}-1][2l])^{1/2}}{[2l+1]}(m)_2^{-11},\\[4mm] 
E_{31}(m)_1&=&\frac{a_1}{a_3}q^{l_{11}-l_{13}} 
\left( 
\frac{[l_{11}-l_{23}]}{[2l+1]}\right)^{1/2}(m)_3^{-11},\\[4mm] 
E_{31}(m)_2&=&\frac{a_2}{a_3}q^{l_{11}-l_{23}} 
\left( 
\frac{[l_{13}-l_{11}]}{[2l]}\right)^{1/2}(m)_3^{-11},\\[4mm] 
E_{31}(m)_3&=&0,\\[6mm] 
E_{13}(m)&=&0,\\[4mm] 
E_{13}(m)_1&=&-a_1q^{l_{23}-l_{11}-1}
\left( 
\frac{[l_{13}-l_{11}]}{[2l+1]}\right)^{1/2}[l_{23}+l_{33}+3]  (m)^{+11},\\[4mm]
E_{13}(m)_2&=&a_2q^{l_{13}-l_{11}-1}
\left(
\frac{[l_{11}-l_{23}]}{[2l]}\right)^{1/2}[l_{13}+l_{33}+3]
 (m)^{+11} \\[4mm]
E_{13}(m)_3&=&a_3
\left\{
\frac{q^{l_{13}-l_{11}-2}}{a_1}
\left(
\frac{[l_{11}-l_{23}+1]}{[2l+1]}\right)^{1/2} [l_{13}+l_{33}+3] (m)_1^{+11}
\right.\\[4mm]
&& \left. +
\frac{q^{l_{23}-l_{11}-2}}{a_2}\left([l_{13}-l_{11}-1][2l]\right)^{1/2}
\frac{[l_{23}+l_{33}+3]}{[2l+1]}(m)_2^{+11}\right\}.
~~~~(3.38) 
\end{eqnarray*} 

A question arising here is  when the representations constructed are irreducible 
and how they are classified. It will be dealt with in the next section.
\\[1cm]
{\bf IV.  TYPICAL AND  NONTYPICAL REPRESENTATIONS OF $U_q[gl(2/1)]$}

  The finite--dimensional representations constructed above are either 
irreducible or indecomposable. We can prove the following proposition. \\[2mm]
{\bf Proposition 4}: {\it The finite--dimensional representations  of  
$U_q[gl(2/1)]$ given in (3.37) and (3.38) are irreducible and called typical 
if and only if the condition 
$$[l_{13}+l_{23}+3][l_{23}+l_{33}+3]\neq 0\eqno (4.1)$$ 
holds.} \\[3mm]
When this condition (4.1) is violated, i.e. one of the following pairs of conditions  
$$[l_{13}+l_{33}+3]= 0  ~~~ \mbox{and} ~~~ [l_{23}+l_{33}+3]\neq 0\eqno (4.2)$$ 
or 
$$[l_{13}+l_{33}+3]\neq 0  ~~~  \mbox{and} ~~~ [l_{23}+l_{33}+3]= 0\eqno (4.3)$$ 
(but not both of them simultaneously) holds, the module $W^{q}$ is no longer irreducible 
but indecomposable. In this case, however, there exists an invariant subspace, say 
$I_k^{q}$, of $W^{q}$ such that the factor representation in the factor module 
$$W_k^{q}:=W^{q}/I_k^{q}\eqno (4.4)$$ 
is irreducible. We call this irreducible representation non-typical in the non-typical 
module $W_k^{q}$. Then, as in \cite {ic96}, it is not difficult for us to prove the 
following assertions.\\ [3mm]
{\bf Proposition 5}:
$$V_3^{q}\subset I_k^{q}\eqno (4.5)$$ 
{\it and} 
$$V_0^{q}\cap I_k^{q}= 0.\eqno (4.6)$$ 
 
From (3.37) -- (4.3) we can easily find all non-typical representations of  
$U_q[gl(2/1)]$ which are classified into two classes.\\[7mm]
{\bf IV.1. {\large Non-typical representations of class 1}} 
 
This class is characterized by the conditions (4.2) which,  
for generic  $q$, take the forms 
$$l_{13}+l_{33}+3= 0\eqno(4.2a) $$ 
and
$$l_{23}+l_{33}+3\neq 0.\eqno (4.2b)$$ 
In the other words, we have to replace everywhere all $m_{33}$ by 
$-m_{13}-1$,  keeping (4.2b) valid. Thus we have the 
following proposition.\\[3mm]  
{\bf Proposition 6}: 
$$ I_1^{q}= V_3^{q}\oplus V_2^{q}.\eqno (4.7)$$ 
Then the class 1 non-typical representations in 
$$W_1^{q}= W_1^{q}([m_{13}, m_{23}, -m_{13}-1])\eqno (4.8)$$ 
are given through (3.31) by keeping the conditions (4.2a) and (4.2b) and 
replacing all vectors belonging to $I_1^{q}$ with 0: 
\begin{eqnarray*} 
~~~~~~~~~~~~~~~~~~~~~~~ E_{23}(m)&=&0,\\[4mm] 
E_{23}(m)_1&=&a_1 
\left( 
\frac{[l_{11}-l_{23}]}{[2l+1]}\right)^{1/2}[l_{23}-l_{13}](m),\\[4mm] 
E_{32}(m)&=&\frac{1}{a_1}\left( 
\frac{[l_{11}-l_{23}]}{[2l+1]}\right)^{1/2}(m)_1,\\[4mm] 
E_{32}(m)_1&=&0.~~~~~~~~~~~~~~~~~~~~~~~~~~~ 
~~~~~~~~~~~~~~~~~~~~~~~~~~~~~~~~~~~~ \label {}{(4.9)} 
\end{eqnarray*} 
\\
{\bf IV.2. {\large Non-typical representations of class 2} }
 
For  this class non-typical representations we must keep the conditions  
$$l_{13}+l_{33}+3\neq 0 \eqno(4.3a)$$
and
$$l_{23}+l_{33}+3=0,\eqno (4.3b)$$ 
derived from (4.3) when the deformation parameter $q$ is 
generic. Equivalently, we have to replace everywhere all $m_{33}$ by 
$-m_{23}$ and keep (4.3a) valid. 
Now the invariant subspace $I_2^{q}$ is determined as follows.\\ [3mm] 
{\bf Proposition 7}:
$$I_2^{q}=V_3^{q}\oplus V_1^{q}.\eqno (4.10)$$ 
 
The class 2 non-typical representations in 
$$W_2^{q}([m_{13},m_{23},-m_{23}])\eqno (4.11)$$ 
are also given through (3.31) but by keeping the conditions (4.3a) and (4.3b) 
valid and replacing all vectors belonging to the invariant subspace $I_2^{q}$ by 0: 
\begin{eqnarray*} 
~~~~~~~~~~~~~~~~~~~~~~ E_{23}(m)&=&0,\\ 
E_{23}(m)_2&=&a_1  
\left(\frac{[l_{13}-l_{11}]}{[2l]}\right)^{1/2}[2l+1](m),\\ 
E_{32}(m)&=&\frac{1}{a_2} 
\frac{([l_{13}-l_{11}][2l])^{1/2}}{[2l+1]}(m)_2,\\ 
E_{32}(m)_2&=&0.~~~~~~~~~~~~~~~~~~~~~~~~~~~~~~~~~~~~~~~~~
~~~~~~~~~~~~~~~~~~~~~~\label{}{(4.12)} 
\end{eqnarray*}

We complete this section with the following statement.\\[3mm] 
{\bf Proposition 8}:
{\it The class of  the finite--dimensional representations determined 
above contains all finite--dimensional irreducible representations of the quantum 
superalgebra $U_q[gl(2/1)]$}.\\[1cm]
{\bf V. CONCLUSION}

     The quantum superalgebra $U_q[gl(2/1)]$ is given as both a Drinfel'd--Jimbo 
deformation of the universal enveloping $U[gl(2/1)]$ and a Hopf superalgebra. 
Using the Hopf algebra structure of  $U_q[gl(2/1)]$ we have constructed all its 
finite (irreducible) dimensional representations in a basis of the even subalgebra 
$U_q[gl(2/1)_0]$. This method combines the advantage of previously suggested 
methods for constructing representations of a classical superalgebra 
\cite{class89} -- \cite{osp32} and a quantum superalgebra  \cite{q94,q95} 
and shows that the method used in the classical case can be 
extended to the quantum deformation case. It proves once again the usefulness 
of knowing a Hopf algebra structure of a quantum group. In particular, using 
a Hopf algebra structure of a quantum superalgebra, $U_q[gl(2/1)]$ in the 
case, one  can calculate in an easier way matrix elements in the induced basis 
and express the latter  in terms of a basis of the even subalgebra. All that could 
not be done via the previously suggested procedure \cite{q94} -- \cite{pq96}.  
Such a description of an induced basis may be physically necessary as in it both 
the origin and the structure of multiplets can be seen clearer. Certainly, the 
method of the present paper can be applied to a bigger quantum superalgebra   
and may be also applicable to the case of multi-parameter deformations, 
for example,  the two-parametric $U_{pq}[gl(2/1)]$. We hope this  method and 
the results obtained here could be useful for physics applications. \\ [1cm]
{\bf Acknowledgement}:  One of the author (N.A.K.) would like to thank  
S. Randjbar--Daemi and the Abdus Salam International Centre for Theoretical 
Physics, Trieste, Italy, for kind hospitality.  This work was supported in part by 
the Vietnam National Research Program for 
Natural Sciences under Grant KT -- 04.1.1.  
  
\end{document}